\documentclass[12pt]{amsart}
\usepackage{amsfonts}
\usepackage{amsmath}
\usepackage{amsxtra}
\usepackage{amssymb,latexsym}
\usepackage[mathcal]{eucal}
\usepackage{amscd}
\input xy
\xyoption{all}

\newcommand{\Bcal}{\mathcal{B}}

\newcommand{\Ecal}{\mathcal{E}}
\newcommand{\Fcal}{\mathcal{F}}

\newcommand{\Pcal}{\mathcal{P}}

\newcommand{\Scal}{\mathcal{S}}
\newcommand{\Tcal}{\mathcal{T}}

\newcommand{\Xcal}{\mathcal{X}}

\newcommand{\ch}{\mathbf{1}}

\newcommand{\Z}{\mathbb{Z}}

\newcommand{\R}{\mathbb{R}}
\newcommand{\C}{\mathbb{C}}
\newcommand{\N}{\mathbb{N}}
\newcommand{\T}{\mathbb{T}}

\newcommand{\XT}{$(X,T)$\ }

\newcommand{\OT}{\mathcal{O}_T}

\newcommand{\OCT}{{\bar{\mathcal{O}}}_T}

\newcommand{\OC}{\bar{\mathcal{O}}}

\newcommand{\al}{\alpha}

\newcommand{\del}{\delta}

\newcommand{\ep}{\epsilon}
\newcommand{\sig}{\sigma}

\newcommand{\La}{\Lambda}

\newcommand{\om}{\omega}
\newcommand{\Om}{\Omega}

\newcommand{\br}{\vspace{3 mm}}
\newcommand{\imp}{\Rightarrow}

\newcommand{\tri}{\bigtriangleup}

\newcommand{\cls}{{\rm{cls\,}}}

\newcommand{\Id}{{\rm{Id}}}

\newcommand{\Homeo}{{\rm{Homeo\,}}}

\swapnumbers \theoremstyle{plain}
\newtheorem{thm}{Theorem}[section]
\newtheorem{cor}[thm]{Corollary}
\newtheorem{lem}[thm]{Lemma}
\newtheorem{prop}[thm]{Proposition}

\theoremstyle{definition}
\newtheorem{defn}[thm]{Definition}

\newtheorem{rmk}[thm]{Remark}

\newtheorem{exe}[thm]{Exercise}

\newtheorem{exas}[thm]{Examples}
\newtheorem{prob}[thm]{Problem}

\numberwithin{equation}{section}

\begin{document}
\title[On two recurrence problems]
{On two recurrence problems}
\author{Michael Boshernitzan and Eli Glasner}

\date{September 28, 2008}


\address{Department of Mathematics\\
     Rice University\\
         Huston, Tx 77005\\
         USA}
\email{michael@rice.edu}

\address{Department of Mathematics\\
     Tel Aviv University\\
         Ramat Aviv\\
         Israel}
\email{glasner@math.tau.ac.il}



\begin{abstract}
We review some aspects of recurrence in 
topological dynamics and focus on two
open problems. The first is an old one concerning
the relation between Poincar\'e and Birkhoff recurrence;
the second, due to Boshernitzan, is about moving recurrence. 
We provide a partial answer to a topological
version of the moving recurrence problem. 
\end{abstract}

\maketitle

\section*{Introduction}
Poincar\'e's recurrence theorem is the first and most basic
theorem of ergodic theory. It asserts that 
given a measure preserving (invertible) dynamical system 
$(X,\mu,\Xcal,T)$ and $A\in \Xcal$ with $\mu(A)>0$,
the set $N(A,A)=\{n\in \Z:\mu(T^n A\cap A)>0\}$
meets every set of the form $(L - L)\setminus\{e\}
=\{n - m: n, m \in L, n \ne m\}$ with infinite $L\subset \Z$.
The proof of this surprising fact is straightforward:
The sets $T^n A, n\in L$, having the same (positive) measure,
can not be all disjoint  (mod $\mu$). 
If $\mu(T^n A \cap T^m A)>0$ then
$\mu(T^{n -m}A\cap A)>0$ hence $n -m \in N(A,A)$.

This basic measure theoretic recurrence theorem has a 
topological counterpart due to 
G. D. Birkhoff.
If $(X,T)$ is a topological dynamical system
($X$ is a compact metric space and 
$T\colon X \to X$ is a homeomorphism of $X$ onto itself),
then there is a {\em recurrent point} in $X$; i.e. there is a point
$x \in X$ such that for every $\ep>0$ there is some
$n \ge 1$ with $d(x, T^nx) < \ep$. A purely topological
proof of this theorem (i.e. one which does not use the
fact that such a system always admits an invariant
probability measure and then applies Poincar\'e's theorem)
follows from the fact that minimal subsystems always exist.
One first applies Zorn's lemma to show that every compact topological system
admits a minimal subset and then 
uses the characterization of a 
point whose orbit closure is minimal as a {\em uniformly recurrent point}
(see Lemma \ref{ur} below). 

Poincar\'e's and Birkhoff's recurrence theorems obtained more
recently a deep and far reaching generalization in the form
of Furstenberg's {\em multiple recurrence theorem}, from which Furstenberg
was able to deduce the famous theorem of Szemer\'edi: 
a subset $A \subset \N$ of positive upper Banach density
contains arbitrarily long arithmetical progressions (see Furstenberg \cite{Fur}). 

In the present work we review some aspects of recurrence in 
topological dynamics as developed by 
Gottschalk and Hedlund \cite{GH},  and Furstenberg \cite{Fur},
including several ``folklore" theorems, and then focus on two
particular open problems. The first is an old one concerning
the relation between Poincar\'e and 
Birkhoff recurrence
(see Problems ($A$), ($A'$) and ($A''$),  or \ref{prob:A},  \ref{prob:A1}
and  \ref{prob:A2}, respectively)
and the second, due to Boshernitzan, is about `moving recurrence'
(see Section \ref{sec:mr}).
While the original ``moving recurrence" problem remains open, we
provide here a partial answer to a topological
version of the problem  (Theorem \ref{mr-thm}).
The paper also contains other new results on topological
recurrence. In particular, in Section \ref{sec:tr} we introduce the 
notion of $r$-Birkhoff sets (approximating that of Birkhoff sets) and 
present some preliminary results concerning these sets.
Finally in Section \ref{Sec-UR} we show that ``absolute moving recurrence"
is equivalent to uniform rigidity.

For related works see \cite{Fur}, \cite{Ru}, \cite{A}, \cite{G}, \cite{W}, 
\cite{AG}, \cite{HY}, \cite{G1} and \cite{BFW}.
A comprehensive review by Frantzikinakis and McCutcheon
on the subject of recurrence in dynamics is forthcoming  \cite{FM}.

This work was done while both authors participated in the
program ``Ergodic Theory and Additive Combinatorics" at  MSRI in the
summer of 2008. We thank MSRI for its support. 
We also thank Benjy Weiss for very helpful discussions.

\br

\section{A reminder of some preliminary definitions and basic results}
Let $(X,T)$ be a dynamical system, where $X$ is a compact
metric space and \mbox{$T\colon X \to X$} is a homeomorphism of $X$ onto itself.
For
$A$ and $B$ subsets of $X$, we let
$$
N(A,B)=\{n \in \Z: T^n A \cap B \ne\emptyset\}.
$$
When $A =\{x\}$ is a singleton we write $N(A,B)=N(x, B)$,
thus
$$
N(x,B)=\{n \in \Z: T^n x \in B\}.
$$
For a point $x \in X$ we write $\OT(x)=\{T^n x : n \in \Z\}$
for the {\em orbit of $x$} and $\OCT(x)$ for the closure of
$\OT(x)$. We say that the system $(X,T)$ is {\em point transitive}
if there is a point $x\in X $ with $\OT(x)$ dense. Such a point
is called {\em transitive}.
We say that the system $(X,T)$ is {\em topologically transitive}
(or just {\em transitive})
if the set $N(U,V)$ is nonempty for every pair $U$ and $V$ of nonempty 
open subsets of $X$. Clearly point transitivity implies
topological transitivity
and using Baire's category theorem one can show that
(for metric systems) conversely, in a topologically transitive system
the set $X_{tr}$ of points whose orbit is dense forms a dense
$G_\del$ subset of 
$X$.
A point $x \in X$ is a {\em recurrent point} if the set $N(x,U)\setminus\{0\}$
is nonempty for every neighborhood $U$ of $X$.  
A dynamical system is called {\em minimal} if every point
is transitive. 

A dynamical system $(X,T)$ is {\em equicontinuous} if the collection of
maps $\{T^n: n \in \Z\}$ is equicontinuous. A minimal equicontinuous
system is called a {\em Kronecker system}.
We have the following classical theorem:

\begin{thm}\label{Kro}
\begin{enumerate}
\item A (metrizable) dynamical system $(X,T)$ is equicontinuous
if and only if there is a compatible metric on $X$ with respect to which 
$T$ is an isometry.
\item
A (metrizable) dynamical system is Kronecker if and only if 
it is isomorphic to a system 
of the form $(G,R_{a})$, where $G$ is a compact second countable
monothetic topological group, $a \in G$ is a topological generator 
(meaning that the cyclic subgroup $\{a^n: n \in \Z\}$ is dense in $G$), 
and the transformation $R_{a}$ is defined by $R_{a}g= ga$.
\end{enumerate}
\end{thm}

There exists a largest monothetic compact topological group
$b\Z$ called the {\em Bohr compactification of the integers}.
If we let  $\phi: \Z \to b\Z$ be the canonical map $\phi: \Z \to b\Z$,
then $a = \phi(1)$ is a topological generator of the group
$b\Z$ and one can associate 
to $b\Z$ a
dynamical system
$(b\Z,R_a)$ as above. This system is minimal and equicontinuous,
but non-metrizable. 

For more information on these basic notions of topological dynamics refer
e.g. to chapter one of \cite{Gl}. 

\br

\section{Some families of subsets of $\Z$ and a famous open problem}
In order to avoid some tedious repetitions we introduce
the notation $\Z_* = \Z \setminus \{0\}$.
\begin{defn}
Let $L \subset \Z_*$.
\begin{enumerate}
\item 
$L$ is a {\em Poincar\'e set} if
whenever $(X,\Xcal,\mu,T)$ is a probability preserving system
and $A \subset X$ is a positive set (i.e. $A \in \Xcal$ and $\mu(A)>0$),
then $N(A,A) \cap L \not=\emptyset$.
Let $\bf{Po}$ denote the collection of Poincar\'e subsets of $\Z_*$.
\item
It is a {\em Birkhoff set} (or {\em a set of topological recurrence})
if whenever $(X,T)$ is a minimal
dynamical system and $U \subset X$ a nonempty open set, then
$N(U,U) \cap L\not=\emptyset$.
Let $\bf{Bir}$ denote the collection of Birkhoff subsets of $\Z_*$.
\item
It is a {\em Bohr set} if whenever $(X,T)$ is a Kronecker
dynamical system and $V \subset X$ a nonempty open set, then
$N(V,V) \cap L\not=\emptyset$.
Let $\bf{Bo}$ denote the collection of Bohr subsets of $\Z_*$.
\end{enumerate}
\end{defn}

\begin{defn}
\begin{enumerate}
\item
A subset ${\Fcal}$ of the power set ${\Pcal}$ of $\Z_*$ is
called a {\em family} when it is hereditary upwards.  That is,
$F_1 \subset F_2$ and $F_1 \in {\Fcal}$ imply $F_2 \in {\Fcal}$.  
\item
If $\Ecal$ is any nonempty subset of $\Pcal$ we let $\Fcal(\Ecal)$
be the smallest family containing $\Ecal$. 
\item
If ${\Ecal}$ is any nonempty subset of $\Pcal$ we let
its {\em dual} $\Ecal^*$ be defined by
$$
\Ecal^* = \{F: F \cap E \neq
\emptyset \ {\text{for all}}\  E \in {\Ecal}\}.
$$
It is easy to check that $\Ecal^*$ is a family and that
$\Fcal(\Ecal)^* = \Ecal^*$.
Clearly ${\Fcal}_1 \subset {\Fcal}_2
\Rightarrow {\Fcal}^*_2 \subset  {\Fcal}^*_1$
and finally for a family ${\Fcal} = \Fcal^{**}$.
\end{enumerate}
\end{defn}

\begin{defn}
\begin{enumerate}
\item
Let $\Ecal_{Po}$ be the collection of all subsets of $\Z_*$
of the form $N(A,A)$, 
whenever $(X,\mu,T)$ is a probability preserving system
and $A \subset X$ is a positive set ($\mu(A)>0$).
We then have $\bf{Po} = \Fcal(\Ecal_{Po})^* = \Ecal_{Po}^* $.
\item
Let $\Ecal_{Bir}$ be the collection of all subsets of $\Z_*$
of the form $N(U,U)$, 
whenever $(X,T)$ is a minimal system
and $U \subset X$ is a nonempty open subset of $X$.
We have $\bf{Bir} = \Fcal(\Ecal_{Bir})^*= \Ecal_{Bir}^*$.
\item
Let $\Ecal_{Bo}$ be the collection of all subsets of $\Z_*$
of the form $N(V,V)$, 
whenever $(X,T)$ is a Kronecker system
and $V \subset X$ is a nonempty open subset of $X$.
We have $\bf{Bo} = \Fcal(\Ecal_{Bo})^*=\Ecal_{Bo}^*$.
\end{enumerate}
\end{defn}

\begin{lem}\label{lem:compare3}
$$
\Ecal_{Po} \supset \Ecal_{Bir} \supset \Ecal_{Bo},
$$
whence
$$
\bf{Po} \subset \bf{Bir} \subset \bf{Bo}.
$$
\end{lem}

\begin{proof}
If $(X,T)$ is a minimal system then the collection
$M_T(X)$ of Borel probability measures on $X$ is never empty
and if $U \subset X$ is open and nonempty,
then $\mu(U) >0$ for every $\mu \in M_T(X)$.
This implies $\Ecal_{Po} \supset \Ecal_{Bir}$.
The inclusion $\Ecal_{Bir} \supset \Ecal_{Bo}$
follows trivially from the definitions.
Finally the last two inclusions follow by duality.
\end{proof}

A beautiful result of Kriz \cite{K} (see also \cite{McC})
shows that $\Fcal(\Ecal_{Po}) \supsetneq \Fcal(\Ecal_{Bir})$.

\begin{prob}[$A$]\label{prob:A}
Is it also true that 
$\Fcal(\Ecal_{Bir}) \supsetneq \Fcal(\Ecal_{Bo})$?
\end{prob}

\begin{rmk}\label{remark}
Since $\bf{Po} = \Fcal(\Ecal_{Po})^*$ and
$\bf{Bir} = \Fcal(\Ecal_{Bir})^*$. Kriz' result is the same as the 
statement $ \bf{Po} \subsetneq  \bf{Bir}$, and Problem ($A$) is equivalent
to the question whether $ \bf{Bir} \subsetneq  \bf{Bo}$.
\end{rmk}


Recall that a collection $\Ecal$ of subsets of $\Z$ is {\em divisible \cite{G-80}}
(or has the {\em Ramsey property} \cite{Fur}) if whenever 
$A$ is in $\Ecal$ and $A = C \cup D$, then at least one of the sets
$C$ and $D$ is in $\Ecal$.

\begin{prop}
\begin{enumerate}
\item
The collection $\Ecal_{Bir}$ forms a filter base; hence $\Fcal(\Ecal_{Bir})$
is a filter. 
\item
The family $\bf{Bir}$ is divisible.
\end{enumerate}
\end{prop}

\begin{proof}
1.\ Let $(X,T)$ and $(Y,T)$ be minimal dynamical systems,
and $U \subset X$, $V\subset Y$ nonempty open sets.
Let $M_0\subset X \times Y$ be a minimal subset of the
product system $(X \times Y, T \times T)$.
Since clearly the $\Z^2$-action defined on $X \times Y$ by the group
$\{T^i \times T^j: (i,j) \in \Z^2\}$ is minimal, there is a pair $(i,j)$
such that $(T^i \times T^j )M_0 \cap  U \times V\ne \emptyset$.
Set $M = (T^i \times T^j )M_0$ and $W = (U\times V) \cap M$.
Then the system $(M, T \times T)$ is minimal, the set $W$
is a nonempty open subset of $M$, and clearly
$$
N(U,U) \cap N(V,V) \supset N(W,W).
$$
Thus $\Ecal_{Bir}$ is indeed a filter base. It follows that $\Fcal(\Ecal_{Bir})$,
which is defined as the smallest family containing $\Ecal_{Bir}$,
is a filter.
We leave it as an exercise to show that the dual family of a 
filter has the Ramsey property 
(and vice versa). Therefore
${\bf{Bir}} = \Fcal(\Ecal_{Bir})^*$ has the Ramsey property.
\end{proof}

\begin{exe}
Show that the families ${\bf {Po}}$ and ${\bf {Bo}}$
are also divisible.
\end{exe}

\br

\section{Examples}\label{sec:exs}

\begin{exas}
\begin{itemize}
\item
For every infinite $L \subset \Z$ the difference set
$ \{n - m : m, n \in L,\ n > m\}$ is a Poincar\'e set. In fact this
statement is just Poincar\'e's recurrence theorem.
\item
Let $p(t)$ be a polynomial with real coefficients taking integer values on
the integers and such that $p(0)=0$. Then the sequence $\{p(n)\}_{n\ge 1}$
is Poincar\'e (see \cite[theorem 3.16]{Fur}). In particular the sequence
$\{n^2\}_{n \ge 1}$ is Poincar\'e. It is easy to see that the sequence
$\{n^2+1\}_{n \ge 1}$ is not Poincar\'e.
\item
Every thick set $L \subset \Z$ is Poincar\'e \cite[page 74]{Fur},
(see definition \ref{syn}.2 below).
For the reader's convenience let us reproduce one of the proofs
given  in \cite{Fur}.  Let $(X,\Xcal,\mu,T)$ be a measure preserving system
and $A \in \Xcal$ with $0 < \mu(A)$. If $A$ is not invariant 
(i.e. $\mu(TA \tri A)>0$) then there exists an $N \ge 1$
such that $\mu(\bigcup_{j=0}^N T^j A) > \mu(\bigcup_{j=0}^\infty T^j A) -
\mu(A)$. Then, for any $M \ge 1$ 
$$
\mu(\bigcup_{j=M}^{M+N} T^j A) > \mu(\bigcup_{j=0}^\infty T^j A) - \mu(A).
$$
This implies that $\mu(\bigcup_{j=M}^{M+N} T^j A \cap A) > 0$, for otherwise
$$
\mu(\bigcup_{j=0}^\infty T^j A) \ge
\mu(A) + \mu(\bigcup_{j=M}^{M+N} T^j A) > \mu(\bigcup_{j=0}^\infty T^j A) .
$$
Thus each sufficiently long interval of integers includes an $n$ with
$\mu(T^n A \cap A)>0$, and this proves that a thick set is Poincar\'e.
\item
Recall that a sequence of integers $\{n_k\}_{k=1}^\infty$ is 
{\em lacunary} if  \mbox{$\inf_k \frac{n_{k+1}}{n_k} > 1$}.
A theorem of Y. Katznelson (see \cite[theorem 5.3]{W})
asserts that a lacunary sequence is never 
Bohr. 
In fact, there is a stronger result  (see  \cite{Pol}, \cite{DeM})
according to which for any lacunary sequences of integers  
{$\{n_k\}_{k=1}^\infty$},  there always exists an irrational  $\alpha\in\R$  
such that  $\inf_k \|n_k\alpha\|>0$.
(We write    $\|x\|=\min_{i\in\Z} |x-i|$\,  for the distance of
a real  $x$  from $\Z$,  the set of integers.) This answered a question of
Erd\"os's in \cite{E}.
\item 
The results in the above example do not extend to slower
growing sequences  (see \cite{AHK}, \cite{Bos}).
An increasing sequence of integers  $\{n_k\}_{k=1}^\infty$  is
called {\em sublacunary} if  \mbox{$\lim_{k\to\infty}\limits \frac{n_{k+1}}{n_k} = 1$}.
There are various results  (\cite{AHK}, \cite{Bos}, \cite{Bou}) which  
indicate  that  for a ``generic"  sublacunary sequence 
$\{n_k\}_{k=1}^\infty$  the limit\, 
$\lim_{k\to\infty}\limits \tfrac1N\sum_{k=1}^N \exp(2\pi i\,n_k\alpha)$\,
exists and vanishes for all real  $\alpha\notin\Z$.
Such sequences are known to be Poincar\'{e}
and, in particular,  Bohr.  (In the above three quoted papers the term
``generic"  has various probabilistic meanings). 

\end{itemize}
\end{exas}

\section{A second formulation of the problem}

\begin{defn}\label{syn}
\begin{enumerate}
\item
A subset $S \subset \Z$ is called {\em syndetic} if there is a
positive integer $N$ such that $S + \{0,1,2,\dots,N\} = \Z$.
\item
A subset $R \subset \Z$ is called {\em thick} (or {\em replete})
if for every positive integer $N$ there is 
an $n \in \Z$ such that
$\{n,n+1,n+2,\dots,n+N\} \subset R$.
\item
A point $x \in X$, where $(X,T)$ is a dynamical system, is
{\em uniformly recurrent} if $N(x,U)$ is syndetic for every
neighborhood $U$ of $x$.
(In \cite{GH} a uniformly recurrent point
is called an {\em almost periodic point}.)
\end{enumerate}
\end{defn}

\begin{exe}
Let $\Scal$ and $\Tcal$ denote the families of
syndetic and thick sets respectively. Show that
$\Scal$ and $\Tcal$ are dual families.
\end{exe}

We have the following important lemma 
(see \cite {GH}). 

\begin{lem}\label{ur}
Let $(X,T)$ be a dynamical system and $x_0 \in X$. Then  
$\OCT(x_0)$ is a minimal subset of $X$ if and only if $x_0$ is
uniformly recurrent.
\end{lem}

\begin{proof}
Suppose first that $x_0 \in X$ has a minimal orbit closure 
$Y=\OCT(x_0)$. Let $U$ be a neighborhood of $x_0$ in $X$.
By minimality there is an $N \ge 1$ such that $Y \subset \bigcup_{j=0}^N T^jU$.
Now given $n \in \Z$ there is some $0\le j \le N$ with $T^n x_0 \in T^jU$. 
Thus $T^{n -j}x_0 \in U$, hence $n-j \in N(x_0,U)$, hence $n = m +j$
for some $m \in N(x_0,U)$. 

Conversely, suppose $x_0$ is uniformly recurrent. Set $Y=\OCT(x_0)$ and
let $M \subset Y$ be a minimal subset of $Y$.
Suppose $M \not = Y$, then $x_0 \not\in M$.
Let $U$ and $V$ be open subsets of $X$ such that 
$x_0 \in U$, $V \supset M$ and $U \cap V=\emptyset$.
Pick some $y_0 \in M$. Then the whole orbit
of $y_0$ is contained in $V$ and for every $N \ge 1$ we can find
$n_N$ with $T^{n_N}x_0$ sufficiently close to $y_0$
to ensure that $T^{n_N}x_0, T^{n_N+1}x_0, \dots, T^{n_N+N}x_0$
are all in $V$. 
This argument shows that the set $N(x_0,U)$ is not syndetic,
contradicting our assumption that $x_0$ is uniformly recurrent.

\end{proof}

\begin{lem}\label{NUU}
Let $(X,T)$ be a dynamical system, $U \subset X$ a
nonempty open subset and $x \in X$. Then
$N(U, U) \supset N(x,U) - N(x,U)$. If moreover, $(X,T)$ is minimal
then $N(U, U)= N(x,U) - N(x,U)$
\end{lem}

\begin{proof}
If $T^mx \in U$ and $T^nx \in U$ then $T^{n-m}T^mx \in U$, so that
$N(U, U) \supset N(x,U) - N(x,U)$. Conversely, if $n \in N(U,U)$, there
is some $y \in U$ with $T^n y \in U$. By minimality there is some 
$m \in \Z$ such that $T^m x$ is sufficiently close to $y$ to ensure
that both $T^m x \in U$ and $T^nT^m x \in U$. Then,
$n = (n+m) - m$ and both $n+m$ and $m$ are in $N(x,U)$. 
\end{proof}

\begin{lem}\label{symbolic}
If $S \subset \Z$ is syndetic then there is a minimal system
$(Y,T)$ and an open nonempty $U\subset Y$ such that 
$S - S \supset N(U,U)$.
\end{lem}

\begin{proof}
Let $\Om=\{0,1\}^\Z$ and $\sig: \Om \to \Om$ the shift
transformation: $ (\sig\om)_n= \om_{n+1}$.
Set $Y' = \OC_{\sig}(\ch_S)$ and $U'=\{\om\in \Om:
\om_0 = 1\}$. It is not hard to check that $Y'$
contains a minimal subset $Y \subset Y'$ such that
$U=Y \cap U'$ is not empty. If $n \in N(U,U)
=\{n: \sig^n U \cap U \ne\emptyset\}$ then
there is a point $y_0 \in U$ with $\sig^n y_0 \in U$.
There exists an $m \in \Z$ such that $\sig^m \ch_S$
is sufficiently close to $y_0$ to ensure that both
$\sig^m \ch_S \in U'$ and $\sig^n\sig^m \ch_S\in U'$.
Thus both $m$ and $n+m$ are in $S$ and $n = n+m - m$
is in $S - S$.
\end{proof}

The next lemma follows easily from the characterization of Kronecker 
systems given in theorem \ref{Kro}; we leave the details to the reader.

\begin{lem}\label{Bo}
Let $(X,T)$ be a Kronecker system and
$V \subset X$ a nonempty open subset. 
Then for every point $x_0\in V$ there exists an open neighborhood 
$x_0 \in V_0 \subset V$ such that 
$$
N(V,V) \supset N(x_0, V) \supset N(V_0,V_0).
$$
Thus, denoting by $\Ecal'_{Bo}$ the collection of subsets of the
form $N(x,V)$, where $(X,T)$ is Kronecker,
$x \in X$ and $V$ is an open neighborhood of $x$, we have
$\Fcal(\Ecal_{Bo})=\Fcal(\Ecal'_{Bo})$, hence
$\bf{Bo} = \Fcal(\Ecal_{Bo})^* = \Ecal^*_{Bo}=  {\Ecal'}_{Bo}^*$.
\end{lem}



\begin{defn}
Let $\al = (\al_1,\al_2,\dots,\al_k)$, a finite sequence of real numbers,
and $\ep >0$ be given. Set
$$
B(\al_1,\al_2,\dots,\al_k; \ep)=\{n \in \Z: \|n\al\| < \ep\}.
$$
Here $\al$ is considered as an element
of the $k$-torus, $\T^k = (\R/\Z)^k$ and for $x \in \R^k$,
$\|x\|$ denotes the Euclidian distance of $x$ from $\Z^k$.
We say that a subset $B$ of $\Z$ is a {\em Bohr neighborhood of zero}
if it contains some $B(\al_1,\al_2,\dots,\al_k; \ep)$.
\end{defn}

Since by Kronecker's theorem the equicontinuous dynamical system
$(\T^k, T)$, where $Tx= x +\al \pmod 1$ with
$\{1,\al_1,\al_2,\dots,\al_k\}$ independent over the rational numbers,
is a minimal system, it follows that every Bohr neighborhood of zero
is in $\Fcal(\Ecal_{Bo})$. (Take $V=B_\ep(0) \subset \T^k$,
so that $B(\al_1,\al_2,\dots,\al_k; \ep)= N(0,V)$.)
With a little more effort one can prove the following
characterizations of Bohr neighborhoods of zero.

\begin{prop}\label{Bo=BN}
The following conditions on a subset $B \subset \Z$
are equivalent:
\begin{enumerate}
\item
$B$ is a Bohr neighborhood of zero.
\item
$B$ is in $\Fcal(\Ecal_{Bo})$; i.e.
$B$ contains a subset of the form $N(V,V)$ where
$(X,T)$ is a Kronecker system and $V$
a nonempty open subset of $X$.
\item
$\cls {\phi(B)}$ is a neighborhood of the zero element in the compact
monothetic group $b\Z$. Here $b\Z$ is the Bohr compactification of the 
integers and $\phi: \Z \to b\Z$ is the natural embedding. 
\end{enumerate}
\end{prop}

\begin{prob}[$A'$]\label{prob:A1}
Given a syndetic subset $S \subset \Z$, is
$S - S$ a Bohr neighborhood of zero? That is, is there
a set $B = B(\al_1,\al_2,\dots,\al_k; \ep)$ with $S - S \supset B$?
\end{prob}

{\bf Claim.}\ 
Problem ($A'$) is a reformulation of Problem ($A$).

\begin{proof}
To see this assume first that the answer to Problem ($A'$) is in the
affirmative.
Let $(X,T)$ be a minimal system and $U \subset X$ a nonempty
open subset. Then, by Lemma \ref{NUU}, $N(U,U) = S - S$, where
$S = N(x_0,U)$ for some (any) $x_0 \in X$. By Lemma \ref{ur},
$S$ is syndetic and by our assumption $S - S$ and therefore also
$N(U,U)$ contain a Bohr neighborhood of zero. By Proposition 
\ref{Bo=BN} we conclude that every $N(U,U)$, i.e. every member
of $\Ecal_{Bir}$, contains a member of $\Ecal_{Bo}$, whence
$\Fcal(\Ecal_{Bir}) = \Fcal(\Ecal_{Bo})$, and $\bf{Bir} = \bf{Bo}$
(see Remark \ref{remark} above).

Conversely, assume now that $\bf{Bir} = \bf{Bo}$,
and let $S \subset \Z$ be a syndetic subset.
By Lemma \ref{symbolic}, $S - S$ contains a set of the form $N(U,U)$ for
some minimal system $(X,T)$ and an open nonempty $U \subset X$.
If the set $S - S$ is not a Bohr neighborhood of zero then,
for every Kronecker system $(Y,T)$ and nonempty open
$V \subset Y$, $N(V,V) \cap N(U,U)^c \ne\emptyset$, and therefore
$N(U,U)^c$ is in $\bf{Bo}$. This contradicts our assumption 
since $N(U,U) \cap N(U,U)^c = \emptyset$ implies that $N(U,U)^c$ is not in
$\bf{Bir}$.
\end{proof}

We do have the following facts:

\begin{thm}
Let  $S \subset \Z$ be a syndetic subset.
\begin{enumerate}
\item (Veech \cite{Veech})
There exists a Bohr neighborhood of zero $B$ such that
$(S - S) \tri B$ is a subset of upper Banach density zero.
\item (Ellis and Keynes \cite{EK})
There exists a Bohr neighborhood of zero $B$
with $S - S + S - s\supset B$ for some $s \in S$.
\end{enumerate}
\end{thm}


Recall that a topological group $G$ is called {\em
minimally almost periodic\/} (MAP) if it admits no nontrivial
continuous homomorphism into a compact group. Or, equivalently,
if it admits no nontrivial minimal equicontinuous action
on a compact space. There are many examples of MAP monothetic
Polish groups (see e.g. \cite{AHK}). 
A topological group $G$ has the {\em fixed point on compacta\/}
property (FPC) if every compact $G$ dynamical system has a fixed point;
see \cite{GrM} and \cite{G}. 
Some authors call this property {\em extreme amenability}.
Recently the theory of Polish groups with 
the fixed point on compacta property received a lot
of attention and new and exciting connections with other branches
of mathematics (like Ramsey theory, Gromov's theory of mm-spaces, and
concentration of measure phenomena) were discovered;
see V. Pestov's book \cite{P}.
In \cite{G} it is shown that the Polish group $G$ of all measurable
functions $f$ from a nonatomic Lebesgue measure space $(\Om,\Bcal,m)$
into the circle $\{z \in \C: |z|=1\}$, with pointwise product and the topology
of convergence in measure, is monothetic and has the
FPC property. 
Of course every topological group with the FPC property is also MAP.
The following problem is posed in \cite{G}.
\begin{prob}[$A''$]\label{prob:A2}
Is there a Polish monothetic group
which is MAP but does not have the fixed point
on compacta property?
\end{prob}
It is shown there that a positive answer to problem $(A'')$
would provide a negative answer to problem $(A')$.

\br

\section{More on topological recurrence}\label{sec:tr}
Let $(X,T)$ be a dynamical system, where $X$ is a compact
metric space and $T : X \to X$ is a homeomorphism of $X$ onto itself.
We fix a compatible metric $d$ on $X$.
Recall the following familiar definition.
A point $x \in X$ is {\em recurrent} if for every 
$\ep>0$
there is a $n\in \Z\setminus \{0\}$ with $d(T^nx,x) <\ep$.
Equivalently, setting  
$$
\phi(x)=\inf \{d(T^nx,x): n\in \Z\setminus \{0\}\},
$$ 
we see that $x$ is recurrent iff $\phi(x)=0$.
More generally, given an {\em infinite} subset $L \subset \Z \setminus\{0\}$, set
$$
\phi_L(x)=\inf \{d(T^nx,x): n\in L\},
$$
and call a point $x \in X$, {\em $L$-recurrent} when $\phi_L(x)=0$.
Let us remark that the role of the metric $d$ in these definitions is not 
essential. It is not hard to show that although the functions $\phi_L$ 
usually depend on the choice of a compatible metric $d$,
the sets of $L$-recurrent points do not. 
We say that a subset $A \subset X$ is
{\em wandering} if there is an infinite set $J \subset \Z$
such that the sets $T^j A; j \in J$, are pairwise disjoint. 
We say that the system $(X,T)$ is {\em non-wandering} if $X$ contains 
no nonempty wandering open subsets. 
Following Furstenberg, \cite[Theorem 1.27]{Fur} , we have:

\begin{thm}\label{five}
\begin{enumerate}
\item
The function $\phi_L$ is upper-semi-continuous. 
\item
The set of $L$-recurrent points is a $G_\del$ subset of $X$.
\item
If $(X,T)$ is non-wandering then the set 
of recurrent points is a dense $G_\del$ subset of $X$.
\item
If there is a $T$-invariant probability measure $\mu$ on $X$
with full support (i.e. $\mu(U)>0$ for every nonempty open $U$)
and $L$ is a Poincar\'e set then the set 
of $L$-recurrent points is a dense $G_\del$ subset of $X$.
\item
If $(X,T)$ is minimal and $L$ is a Birkhoff set then the set 
of $L$-recurrent points is a dense $G_\del$ subset of $X$.
\end{enumerate}
\end{thm}

\begin{proof}
We leave the proofs of the claims (1) and (2) as an exercise.
For (3) see Furstenberg, \cite[Theorem 1.27]{Fur} (or adapt
the
following proof).
For the proof of claim (4) we first recall that an upper-semicontinuous function
on $X$ has a dense $G_\del$ set of continuity points.
Let $X_L \subset X$ be the dense $G_\del$ set of continuity
points of $\phi_L$. Suppose $\phi_L(x_0)=a >0$ for some $x_0\in X_L$. 
Then, by continuity, there is a $0< \delta < a/4$ such that $\phi(x) > a/2$ 
for every $x$ is an open ball $U$ of radius $\del$ around $x_0$.
Since $\mu(U)>0$ and $L$ is Poincar\'e we have 
$L \cap N(U,U) \ne \emptyset$.
For $n$ in this intersection there are $u_1,u_2 \in U$ with
$T^{n} u_1 = u_2$, hence $d(u_1,T^{n} u_1)< a/2$.  
In particular $\phi_L(u_1) < a/2$. This contradicts our 
choice of $U$ and we conclude that $\phi_L(x)=0$
for every $x \in X_L$. This completes the proof of claim (4).
A similar argument will prove claim (5).
\end{proof}


In the next two theorems we establish several characterizations of
Birkhoff sets. We will use the following lemma which is valid for every minimal system.

\begin{lem}\label{M}
Let $(X,T)$ be a minimal system and let $\eta>0$ be given.
Then there exists a positive integer $M \ge 1$ such that
for every $x \in X$ the set $\{T^j x \}_{j=0}^M$ is $\eta$-dense in $X$;
i.e. for every $x'\in X$ there is some $0 \le j \le M$ with 
$d(x',T^jx)< \eta$.
\end{lem}

\begin{proof}
Assuming the contrary we would have for each $n$,
points $x_n, y_n \in X$ such that $B_\eta(y_n)\cap
\{T^jx_n\}_{j=0}^n =\emptyset$.
By compactness there are convergent subsequences 
say, $x_{n_j} \to x$ and  $y_{n_j} \to y$.
By minimality there is a positive $m \ge 1$ such that
$d(T^m x, y) < \eta/3$. We now choose $j$ so large that:
(i) $n_j > m$, (ii) $d(y,y_{n_j}) < \eta/3$, and
(iii) $x_{n_j}$ is sufficiently close to $x$ to ensure that
$d(T^m x_{n_j}, T^m x) < \eta/3$.
With this choice of $j$ we now have:
$$
d(T^m x_{n_j}, y_{n_j}) < d(T^m x_{n_j}, T^m x) + d(T^m x, y)
+ d(y, y_{n_j}) < \eta/3+\eta/3+\eta/3=\eta.
$$
Since $n_j > m$, this contradicts the choice of $x_{n_j}$ and
$y_{n_j}$.
\end{proof}

\begin{thm}\label{rec-thm}
The following conditions on a subset $L \subset \Z_*$ are equivalent.
\begin{enumerate}
\item
$L$ is Birkhoff.
\item
$L \cap (S - S) \ne \emptyset$ for every syndetic subset $S \subset \Z$.
\item
For every minimal dynamical system $(X,T)$, the set of
$L$-recurrent points is dense and $G_\del$.
\item
For every dynamical system $(X,T)$ and $\ep>0$
there are $x \in X$ and $m \in L$ with $d(T^mx,x) < \ep$.
\end{enumerate}
\end{thm}

\begin{proof}
From Lemmas \ref{NUU} and \ref{symbolic} we easily deduce
the equivalence of properties (1) and (2).
The implication (1) $\imp$ (3) is proven in Theorem \ref{five}.5.
Next assume (3).  Given a minimal system $(X,T)$ and
a nonempty open subset $U \subset X$ we clearly have
$L \cap N(U,U)\ne\emptyset$, whence $L$ is Birkhoff.
Thus we have (3) $\imp$ (1).

As every dynamical system has a minimal subsystem we clearly have  
(3) $\imp$ (4). Finally we show that (4) implies (3).
Let $(X,T)$ be a minimal system. For $\ep>0$ set
$$
V_L(\ep)=\{x \in X: \exists\ m \in L \ {\text{with}}\ d(T^mx,x)< \ep\}.
$$
Clearly $V_L(\ep)$ is open and assuming (4) we know that it is nonempty.
Given $\eta >0$ there is, by Lemma \ref{M}, $M \ge 1$ such that for every
$x \in X$ the set $\{T^i x\}_{i=0}^M$ is $\eta$-dense. 
Let $\del >0$ be such that $d(x,x') < \del$ implies $d(T^ix,T^ix')< \ep$
for every $0 \le i \le M$. It now follows that for every $x \in V_L(\del)$,
we have  $\{T^i x\}_{i=0}^M\subset V_L(\ep)$, and consequently, that
$V_L(\ep)$ is $\eta$-dense. Since $\eta$ is arbitrary, we conclude that
$V_L(\ep)$ is dense. By Baire's theorem we conclude that
$X_0 = \bigcap_{\ep>0} V_L(\ep)$ is a dense $G_\del$ subset of $X$.
Clearly every $x \in X_0$ is $L$-recurrent.
\end{proof}

In order to achieve additional characterizations for Birkhoff sets we introduce 
the following definition. For $r\in\N$ we denote $\N_r=\{1,2,\ldots,r\}$.


\begin{defn}\label{def:color}
Let  $r\in \N$.  A subset  $L\subset\Z_*$  is said to be  {\em $r$-Birkhoff}
(notation:  $L\in{\bf Bir}_r$)
if the following two equivalent conditions hold:
\begin{enumerate}
\item 
For every sequence  $\{z_i\}_{i\in\Z}$  over  $\N_r$,
there are  $m\in L$ and $i\in\Z$ such that $z_i=z_{i+m}$.
\item  
For every coloring   $c\colon \Z\to\N_r$  there are
$i,j\in \Z$,  with  $c(i)=c(j)$  and $i-j\in L$. 
\end{enumerate}
\end{defn}

\begin{rmk}
In the above definition one can replace $\Z$ by $\N$.
\end{rmk}
\begin{thm}\label{rec-thm2}
The following conditions on a subset $L \subset \Z_*$ are equivalent:
\begin{enumerate}
\item
$L$ is Birkhoff.
\item
For any compact metric space $Z$, every sequence $\{z_i\}_{i \in \Z}$,
with $z_i \in Z$,
and every $\ep >0$, there are $m \in L$ and $i \in \Z$ such that 
$d(z_i, z_{i +m}) < \ep$.
\item
$L$  is $r$-Birkhoff\, for all\,  $r\in\N$.
\end{enumerate}
\end{thm}

Thus, by the above theorem,   {\bf Bir} $=\bigcap_{r\in\N}\, ${\bf Bir}$_r$. 

\begin{proof}
(1) $\imp$ (2): Suppose $L$ is Birkhoff and let $Y = \OCT(\zeta)$, where
$\Om=Z^\Z$, $T:\Om \to \Om$ is the shift and the
element $\zeta\in \Om$ is defined by $\zeta(i)=z_i$.
Let $M \subset Z$ be a minimal subset. Applying Theorem \ref{rec-thm}
we see that there is a point $x \in M \subset Y$ which is $L$-recurrent.
Fix a compatible metric $d$ on $\Om$ and let $0 < \del$ be
such that $d(\om,\om')< \del$ implies $d(\om(0),\om'(0))<\ep$.
Let $m \in L$ be such that $d(T^m x,x) < \del$. Let  $i \in \Z$
be chosen so that $T^i\zeta$ is sufficiently close to $x$ to ensure
that also $d(T^m T^i \zeta, T^i\zeta)< \del$. By our choice of
$\del$ we have $d(z_{m+i},z_i)=d(\zeta(m + i),\zeta(i))< \ep$.

(2) $\imp$ (3):  Take $Z =\N_r=\{1,2,\dots,r\}$. Let $d(i,j)=\del_{ij}$
for $i,j\in \N_r$  and take $\ep =1/2$. Then
$d(z_i, z_{i +m}) < \ep$ implies $z_i = z_{i +m}$.

(3) $\imp$ (1): We will show that condition (4) in Theorem \ref{rec-thm}
is satisfied. So let $(X,T)$, a minimal system, and $\ep >0$ be given.
Let $\{V_i\}_{i=1}^r$ be an open cover of $X$ by balls of radius
$\ep/2$. Fix $x_0 \in X$ and choose a sequence $z_i \in \N_r$
such that $T^i x_0 \in U_{z_i}$ for every $i \in \Z$.
By (3) we have $m \in L$ and $i\in \Z$ such that $z_{m+i}=z_i=j$,
whence $T^i x_0$ and $T^{m+i} x_0$ are both in $U_j$.
Thus $d(T^m T^i x_0, T^i x_0) < \ep$, and taking $x = T^i x_0$
we have the required $x$.
\end{proof}

\begin{rmk}
The last condition in Theorem \ref{rec-thm2} can be formulated as
a coloring property:
For every $r$ and every coloring $c: \Z \to \{1,2,\dots,r\}$ there are 
$i, j \in \Z$, with $c(i) = c(j)$ and $i -j \in L$. 
See \cite{W} for a graph theoretical interpretation of this coloring property.
\end{rmk}

We will now consider some basic properties of $r$-Birkhoff sets.
The first statement we leave as an easy exercise.

\begin{exe}[A Compactness Principle]
For  $r\geq1$,  every  $r$-Birkhoff set  contains
a finite $r$-Birkhoff subset.
\end{exe}

For any  $r\in \N$,  each of the sets\,  $k\N_r=\{k,2k,\ldots,rk\}$, 
$k\in\N$,  is  \text{$r$-Birkhoff}.  
Indeed,  let  
$(z_i)$  be an arbitrary sequence in   $\N_r$.  
Since  card$(k\N_{r+1})=r+1>r
=$ card$(\N_r)$, 
there are   $i, j\in k\N_{r+1}$, $i\neq j$,  such that  $z_i=z_j$.  
Assuming, with no loss of generality,  that   $m=j-i>0$,  we get  $z_i=z_{i+m}$,  
with some  $m\in k\N_r$, completing the proof  (see
Definition \ref{def:color}, first condition).  

On the other hand, for finite subsets  $M\subset \Z_*$  the  following
implication holds:
\[
\text{card}(M)=r\geq1\quad \implies \quad M\notin \text{\bf Bir}_{r+1}.
\]
With no loss of generality we may assume that  $M\subset \N$  (by replacing
$M$  by the set   
$(M\cup(-M))\cap\N$).
Construct a sequence  $\{z_i\}_{i\in\Z}$ over the set  
$\N_{r+1}=\{1,2,\ldots,r+1\}$ as follows.
For $i\leq0$,  set  $z_i=1$;  for $i\geq1$,  set inductively: 
$$
 z_i=\min X_i, \quad \text{where }\ X_i=\{x\in\N_{r+1}\mid  x\neq z_{i-m}, \text{ for all } m\in M\}.
 $$
(Clearly,  $X_i\neq\emptyset$ for $i\geq1$, because  card$(M)=r<r+1=$\,card$(\N_{r+1})$).
The above construction implies that 
 $z_i=z_{i+m}$  has no solutions  in  $i\in\N$  and  $m\in M$.
It follows that  $M\notin$\ {\bf Bir}$_{r+1}$  (see  Definition  \ref{def:color} 
and the subsequent remark).
 
We conclude that the sets   $k\N_r=\{k,2k,\ldots,rk\}$, 
$k,r\in\N$,  provide examples of  \mbox{$r$-Birkhoff} sets which are not   $(r+1)$-Birkhoff.
More refined examples will be provided 
next
(see \eqref{eq:lr} below).

\begin{defn}
A subset  $M\subset Z_*$  is called {\em stably} $r$-Birkhoff  
(notation:  $M\in \text{\bf Bir}'_r$)  if for 
every finite subset  $F\subset \Z$,
the difference set  $M\setminus F$  is  $r$-Birkhoff.
\end{defn}
Define the sets   
\begin{equation}\label{eq:lr}
L_r=\Big\{n(r+2)^k\mid n\in\{1,2,\ldots,r\}, k\geq0\Big\}\subset\N, \quad\text{for }\ r\in\N.
\end{equation}

We claim that, for every  
$r \ge 2$,
the set $L_r$ 
\begin{enumerate}
\item[]
\begin{enumerate}
\item is lacunary;
\item is stably  $r$-Birkhoff;
\item is not $(r+1)$-Birkhoff.
\end{enumerate}
\end{enumerate}

The fact that  $L_r$  is lacunary is 
clear.  In fact, if  $\{x_1<x_2<\ldots\}$    is the
linear ordering of  $L_r$, then\,   
$\min_{k\geq1} \frac{x_{k+1}}{x_k}=\frac r{r-1}$, for  $r\geq2$.

The set   $L_r$  is stably $r$-Birkhoff  because  $L_r$  can be represented as a disjoint infinite union   $L_r=\bigcup_{k\geq0} L_{r,k}$  where each\,  $L_{r,k}=(r+2)^k \N_r$\,  is   $r$-Birkhoff  (as proved earlier).

Finally, to prove that $L_r$  is not $(r+1)$-Birkhoff, 
define a sequence  $\{z_k\}_{k\in\Z}$
over  $\N_{r+1}$  by the condition   $z_i\equiv i$ (mod$(r+1)$).  
We claim  that   $z_i=z_{i+m}$  has no solution
in   $m\in L_r$  and  $i\in \Z$.  Indeed, otherwise   
$$
i\equiv i+m\text{ (mod }(r+1))  \implies  m\equiv 0\text{ (mod }(r+1)) \implies  
\tfrac m{r+1}\in\Z
$$
which is impossible  (see \eqref{eq:lr}). This completes the proof that    $L_r\notin$\ {\bf Bir}$_r$
(see the first condition in Definition \ref{def:color}).

The fact 
that
the sets  $L_r\in\text{\bf Bir}'_r$  are lacunary should be compared with 
the fact that no set  in 
$\text{\bf Bo}$ (which by Lemma \ref{lem:compare3} contains
$\text{\bf Bir}$) is lacunary
(see the examples in Section \ref{sec:exs}).

\section{The moving recurrence problem}\label{sec:mr}

The following question was recently posed by Boshernitzan,
and is still open. 

\begin{prob}[$B$]
Let $(X,T)$ be a dynamical system, $\mu \in M_T(X)$ a $T$-invariant
probability measure on $X$ and $(n_k)$ an infinite sequence of 
nonzero integers. 
Define
$$
\psi_{(n_k)}(x) = \inf_{k\ge 1} d(T^{n_k}x, T^{n_k +k}x).
$$
Is it true that $\psi{(n_k)}(x)=0$, $\mu$-a.e.?
\end{prob}

In this section we prove a topological analogue using the tools
developed in the previous sections.

\begin{defn}
For a sequence $(n_k)$ of elements of $\Z$, let
$$
\psi_{(n_k)}(x) = \inf_{k\ge 1} d(T^{n_k+k} x,T^{n_k} x).
$$
More generally,
given two sequences $(n_k)$ and $(r_k)$ of elements of $\Z$, let
$$
\psi_{(n_k,r_k)}(x) = \inf_{k\ge 1} d(T^{n_k+r_k} x,T^{n_k} x).
$$
We say that a point $x\in X$ is {\em $(n_k)$-moving recurrent} if
$\psi_{(n_k)}(x) =0$. It is {\em $(n_k,r_k)$-moving recurrent} when
$\psi_{(n_k,r_k)}(x) =0$.
Note that $\psi_{(n_k)}=\psi_{(n_k,k)}$.
\end{defn}

Again we have:
\begin{lem}\label{usc}
The function $\psi_{(n_k,r_k)}$ is upper-semi-continuous and the set
of $(n_k,r_k)$-moving recurrent points is a $G_\del$ subset of $X$.
\end{lem}





\begin{thm}\label{mr-thm}
Let $(r_k)$ be a Birkhoff set.
Then for every sequence $(n_k)$,
and every minimal dynamical system $(X,T)$ the set of
$(n_k,r_k)$-moving recurrent points is dense and $G_\del$.
In particular, taking $r_k =k$ we see that for every minimal 
dynamical system $(X,T)$ the set of
$(n_k)$-moving recurrent points is dense and $G_\del$.
\end{thm}

\begin{proof}
{\bf Step 1:} Let $X_0\subset X$ denote the dense $G_\del$ set of 
continuity points of $\psi_{(n_k,r_k)}$ (Lemma \ref{usc}).
Let $x_0 \in X_0$ and assume that $\psi_{(n_k,r_k)}(x_0)=2\ep > 0$.
Since $x_0$ is a continuity point we can find a ball $U$ 
around $x_0$ such that $\psi_{(n_k,r_k)}(x)> \ep$ for every 
$x \in U$.

\br

{\bf Step 2:} 
We will show that the set
$$
V(\ep)=\{x \in X: \psi_{(n_k,r_k)}(x) < \ep\}
$$
is dense.

Fix $\eta > 0$ and use Lemma \ref{M} to find $M \ge 1$ such that 
for every $x \in X$ the set $\{T^j x \}_{j=0}^M$ is $\eta$-dense in $X$.
Next choose  $\del > 0$ such that $d(x,x') < \del$ implies 
$d(T^jx, T^jx') < \ep$ for every $0 \le j \le M$.

Next observe that there exit $x \in X$ 
and $k \ge 1$ with $d(T^{n_k + r_k}x,T^{n_k}x) < \del$.
In fact, since $(r_k)$ is Birkhoff,
we can (applying Theorem \ref{rec-thm}) pick an $(r_k)$-recurrent point
$x' \in X$ and then find $k \ge 1$ with $d(T^{r_k} x', x') < \del$.
Set $x = T^{-r_k}x'$, so that $x' = T^{r_k}x$ and
$$
d(T^{n_k + r_k}x,T^{n_k}x) = d(T^{r_k}x', x') < \del.
$$

Now, by the choice of $\del$, we have:
$$
d(T^{n_k + r_k}T^j x,T^{n_k}T^j x) < \ep,
\ {\text{for all}} \ 0 \le j \le M.
$$
Thus $\psi_{(n_k,r_k)}(T^j x) < \ep$  for all  $0 \le j \le M$
and we conclude that $V(\ep)$ is $\eta$-dense. As this holds
for every $\eta> 0$ we conclude that $V(\ep)$ is dense.

\br

{\bf Step 3:} 
We now have $U \cap V(\ep) \ne \emptyset$ and we achieved  
the conflict $\ep < \psi_{(n_k,r_k)}(x) < \ep$ for any point
$x$ in this intersection. 
Since the assumption $\psi_{(n_k,r_k)}(x_0) > 0$ leads to a conflict
we conclude that $\psi_{(n_k,r_k)}(x_0) = 0$ for every $x_0 \in X_0$,
as required.
\end{proof}

Recall the following definition from \cite{GW}.

\begin{defn}
A dynamical system $(X,T)$ is called an {\em $M$-system}
if (i) it is topologically transitive and (ii) the union of the minimal
subsystems of $X$ is dense in~$X$.
\end{defn}

The class of $M$-systems is very large, e.g. it contains every
topologically transitive system with a dense set of periodic points.
(The latter systems are called {\em chaotic in the sense of Devaney},
or {\em $P$-systems}.)

\begin{cor}
Let $(X,T)$ be an $M$-system and $(n_k)$ an infinite
sequence in $\Z$. Then there is a dense $G_\del$ subset 
$X_0 \subset X$ such that $\psi_{(n_k)}(x) = 0$ for every
$x \in X_0$. In particular the set $X_{tr} \cap X_0$ of
$(n_k)$-moving recurrent transitive points is a dense $G_\del$
subset of $X$.
\end{cor}

\begin{proof}
The set $X_0 = \{x \in X : \psi_{(n_k)}(x)=0\}$ is a $G_\del$ subset of $X$ . 
By Theorem \ref{mr-thm} for every minimal subset $M \subset X$ the set
$M_0 = M \cap X_0$ is a dense $G_\del$ subset of $M$. Thus
$\bigcup\{M_0: M \ {\text{is a minimal subset of }}\  X\}\subset X_0$ is dense in
$\bigcup\{M: M \ {\text{is a minimal subset of }}\ X\}$.
In turn, the latter is dense in $X$ and it follows that $X_0$ 
is dense in $X$. Finally, as the set $X_{tr}$ of transitive points in an
$M$-system is always dense and $G_\del$ we conclude that so is
$X_{tr} \cap X_0$.
\end{proof}

\br

\section{Absolute moving recurrence}\label{Sec-UR}

\begin{defn}
We will say that a system $(X,T)$ is {\em absolutely moving
recurrent} if for every infinite sequence $(n_k)\subset \Z$,
$\psi_{(n_k)} \equiv 0$ (i.e. every point of $X$ is $(n_k)$-moving recurrent). 
\end{defn}

Recall the following definition from \cite{GM}:

\begin{defn}
A dynamical system $(X,T)$ is called {\em uniformly rigid}
if there exists a sequence $m_i \nearrow \infty$ in $\Z$
such that
$$
\lim_{i \to \infty} \sup_{x \in X} d(x, T^{m_i}x) =0.
$$
\end{defn}

For any dynamical system $(X,T)$ let
$$
\La(X,T)={\rm {unif{\text{-}}}} \cls \{T^n: n\in \Z\} \subset
\Homeo (X)
$$
be the uniform closure of the powers of $T$ in the Polish group
$\Homeo (X)$. Of course $\La$ is a Polish monothetic group,
and the system \XT is  uniformly rigid iff it is not discrete.

\begin{thm}
A 
topologically transitive
dynamical system $(X,T)$ is absolutely moving recurrent
if and only if it is uniformly rigid.
\end{thm}

\begin{proof}
Suppose first that $(X,T)$ is uniformly 
rigid with
$T^{m_i} \overset{\text{unif}}{\longrightarrow} \Id$. Let $(n_k)\subset \Z$
be an arbitrary infinite sequence. Then for every $\ep>0$ there exists an
$i_0$ such that for $i > i_0$, $d(T^{m_i}x,x) < \ep$ for every $x \in X$.
In particular then,
$$
d(T^{n_{m_i}+ m_i}x,T^{n_{m_i}}x) = 
d(T^{m_i}(T^{n_{m_i}}x),T^{n_{m_i}}x)< \ep.
$$
Thus, $\liminf_{k}d(T^{n_k+ k}x,T^{n_k}x)=0$ for every $x \in X$,
and we have shown that $(X,T)$ is absolutely moving recurrent.
(Note that 
topological transitivity
is not needed in this direction.)

Conversely, suppose $(X,T)$ is not uniformly 
rigid, that is:
$$
{\text{there exists}}\  \ep_0 > 0, \ {\text{such that}}\  \forall k, \exists x_k \in X,
\ {\text{with}} \ d(T^kx_k,x_k) > \ep_0.
$$
Fix $x_0 \in X_{tr}$ and for each $k$ choose $n_k\in\N$ such that
$T^{n_k}x_0$ is sufficiently close to $x_k$ to ensure
that also $d(T^{n_k+k}x_0,T^{n_k}x_0) = 
d(T^kT^{n_k}x_0,T^{n_k}x_0)> \ep_0$.
For the sequence $(n_k)$ we have $\psi_{(n_k)}(x_0) \ge \ep_0$,
hence $(X,T)$ is not absolutely moving recurrent.
\end{proof}

\br


\end{document}